\documentclass{article}
\usepackage{amsmath,amsthm,amssymb,authblk}
\usepackage[mathscr]{eucal}

\newtheorem{theorem}{Theorem}[section]

\newtheorem{proposition}[theorem]{Proposition}
\newtheorem{lemma}[theorem]{Lemma}

\newtheorem{remark}[theorem]{Remark}
\numberwithin{equation}{section}

\begin{document}
\title{Permutation-like Matrix Groups with a Maximal Cycle of Prime Square Length}
\author{Guodong Deng,\quad Yun Fan\\
\small School of Mathematics and Statistics\\
\small Central China Normal University, Wuhan 430079, China}
\date{}
\maketitle

\insert\footins{\noindent\footnotesize{\it Email address}:
a2b3c4d5deng@163.com (Guodong Deng), yfan@mail.ccnu.edu.cn (Yun Fan).}

\begin{abstract}
A matrix group is said to be permutation-like
if any matrix of the group is similar to a permutation matrix.
G. Cigler proved that, if a  permutation-like matrix group
contains a normal cyclic subgroup which is generated by a maximal cycle
and the matrix dimension is a prime, then the group is similar to
a permutation matrix group.
This paper extends the result to the case
where the matrix dimension is a square of a prime.

\medskip{\em Key words}:~ matrix group, matrix similarity,
permutation-like group, permutation matrix group.

\medskip{\em Mathematics Subject Classification 2010}:~ 15A18, 15A30, 20H20.
\end{abstract}

\section{Introduction}
A multiplicative group consisting of complex invertible matrices of
size $n\times n$ is said to be a matrix group of dimension $n$.
A matrix group ${\cal G}$ is said to be permutation-like
if any matrix of ${\cal G}$ is similar to a permutation matrix,
see \cite{C05, C07}.
If there exists an invertible matrix $Q$ such that
$Q^{-1}AQ$ is a permutation matrix for all $A\in{\cal G}$,
then we say that ${\cal G}$ is similar to a permutation matrix group,
or ${\cal G}$ is a permutation matrix group for short.
A matrix is called a maximal cycle if it is similar to a
permutation matrix corresponding to a
cycle permutation with cycle length equal to the dimension.
G. Cigler in \cite{C07} showed that a permutation-like matrix group
is not a permutation matrix group in general, and suggested
a conjecture as follows.

\noindent{\bf Conjecture.}~
{\em A permutation-like matrix group containing a maximal cycle
is similar to a permutation matrix group.}

\noindent
G. Cigler in \cite{C07} proved it affirmatively in two cases:
the dimension $\le 5$,
or the dimension is a prime integer and the cyclic subgroup
generated by the maximal cycle is normal.

In this paper we extend the result of \cite{C07} to the case where the
length of the maximal cycle is a square of a prime.

\begin{theorem}\label{theorem}
Let ${\cal G}$ be a permutation-like matrix group of dimension $p^2$
where $p$ is a prime.
If ${\cal G}$ contains a maximal cycle $C$ such that
the subgroup $\langle C\rangle$ generated by $C$ is normal in ${\cal G}$,
then ${\cal G}$ is a permutation matrix group.
\end{theorem}

In Section 2 we state some preliminaries as a preparation.
The theorem will be proved in Section~3.

\section{Preparation}
The complex field is denoted by ${\Bbb C}$. For a positive integer $n$,
by ${\Bbb Z}_n^*$ we denote the multiplicative group consisting of
units of the residue ring ${\Bbb Z}_n$ of the integer ring
${\Bbb Z}$ modulo $n$.
A diagonal blocked matrix $\begin{pmatrix}B_1\\ & \ddots\\ &&B_k\end{pmatrix}$
is denoted by $B_1\oplus\cdots\oplus B_k$ for short.
The identity matrix of dimension $n$ is denoted by $I_{n\times n}$.
All complex invertible matrices of dimension $n$ consist the so-called
general linear group, denoted by ${\rm GL}_n({\Bbb C})$.
We denote the characteristic polynomial of a
complex matrix $A$ by ${\rm char}_A(x)$.

\begin{lemma}\label{permutable}
The following two are equivalent to each other:

(i)~  $A$ is similar to a permutation matrix;

(ii)~ $A$ is diagonalizable and ${\rm char}_A(x)=\prod_{i}(x^{\ell_i}-1)$.

\noindent
If it is the case, then each factor $x^{\ell_i}-1$ of ${\rm char}_A(x)$
corresponds to exactly one $\ell_i$-cycle of the cycle decomposition
of the permutation of the permutation matrix.
\end{lemma}

{\bf Proof.}~ It is clear. \qed

\medskip
We'll apply the lemma to the case where $A^{p^2}=1$ with $p$ being a prime,
at that case it is easy to check the condition (ii) of the lemma,
because there are only three divisors $1,p,p^2$ of~$p^2$
which form a chain with respect to the division relation,
and $x^{p^2}-1=\Phi_1(x)\Phi_p(x)\Phi_{p^2}(x)$ where
$\Phi_{p^i}(x)$ denotes the $p^i$'th cyclotomic polynomial.

\medskip
Let $C\in{\rm GL}_n({\Bbb C})$ be a maximal cycle of dimension $n$,
and $\lambda$ be a primitive $n$'th root of unity.
Then there is a basis of the vector space ${\Bbb C}^n$:
\begin{equation}\label{diag basis}
 e_0,~e_1,~\cdots,~e_{n-1},
\end{equation}
such that $C^ie_j=\lambda^{ij}e_j$ for all $i, j=0,1,\cdots,n-1$;
with respect to this basis, $C$ is a diagonal matrix
$C=1\oplus\lambda\oplus\cdots\oplus\lambda^{n-1}$;  see \cite[\S4]{C07}.
And, such basis is unique up to non-zero scales,
since the $1$-dimensional subspace ${\Bbb C}e_j$, for $j=0,\cdots,n-1$,
is just the eigen-subspace of the eigenvalue $\lambda^j$,
for $j=0,\cdots,n-1$ respectively, of the matrix $C$;
or, in representation theoretic notations, ${\Bbb C}e_j$ for $j=0,\cdots,n-1$ are
just all irreducible modules of the cyclic group
$\langle C\rangle$, see \cite[\S15 Example 1]{AB}.

Taking any non-zero complexes $c_0,\cdots,c_{n-1}$ and
setting $f=\sum_{j=0}^{n-1} c_je_j$, we obtain another basis of ${\Bbb C}^n$:
\begin{equation}\label{permutation basis}
 f, ~ Cf,~\cdots,~C^{n-1}f,
\end{equation}
with which $C$ is a cycle permutation matrix, see \cite[Lemma 4.1]{C07}.

Let $B\in{\rm GL}_n({\Bbb C})$ with ${\rm ord}(B)=\ell$,
where ${\rm ord}(B)$ denotes the order of $B$.
Assume that $B$ normalizes the group~$\langle C\rangle$. 
Since the automorphism group
of the cyclic group $\langle C\rangle$ is isomorphic to ${\Bbb Z}_n^*$,
there is an $r\in{\Bbb Z}_n^*$ such that
\begin{equation}\label{conjugation}
 B^{-1}C^iB=C^{ri},\qquad \forall~~ i\in{\Bbb Z}_n\,;
\end{equation}
thus the action by conjugation of $B$ on $\langle C\rangle$
is determined by the action of $\mu_r$ on ${\Bbb Z}_n$,
where $\mu_r(a)=ra$ for all $a\in{\Bbb Z}_n$. Further,
$$
C^iBe_j=B\cdot B^{-1}C^iBe_j=B\cdot C^{ri}e_j=B\cdot \lambda^{rij}e_j
=\lambda^{rij}Be_j;
$$
taking $i=1$, we see that 
$Be_j$ is an eigenvector of the eigenvalue $\lambda^{rj}$ of $C$, i.e.
\begin{equation}\label{permute eigens}
 B{\Bbb C}e_j={\Bbb C}e_{rj},\qquad j=0,1,\cdots,n-1.
\end{equation}
Thus, $B$ permutes the eigen-subspaces
$\big\{{\Bbb C}e_0, {\Bbb C}e_1, \cdots, {\Bbb C}e_{n-1}\big\}$ of $C$
in the same way as $\mu_r$ permutes ${\Bbb Z}_n$.

Let $\Gamma_1,\Gamma_2,\cdots,\Gamma_m$ be orbits of
the action of the group $\langle B\rangle$
on the set of $1$-dimensional subspaces
$\{{\Bbb C}e_1,{\Bbb C}e_2,\cdots,{\Bbb C}e_n\}$.
Assume that the length of $\Gamma_k$ is~$n_k$ for $k=1,\cdots,m$.
Since $|\langle B\rangle|=\ell$ where $|\langle B\rangle|$
denotes the order of the group $\langle B\rangle$,
we have that $n_k|\ell$ for $k=1,\cdots,m$.

\begin{lemma}\label{B block} Let notation be as above.

(1)~ For each $k$, take any one ${\Bbb C}e_{j_k}\in\Gamma_k$ and
any non-zero $e'_k\in{\Bbb C}e_{j_k}$,
set $V_k=\bigoplus_{h=0}^{n_k-1}B^h{\Bbb C}e'_{k}$ and
${\cal E}_k=\{e'_k,Be'_k,\cdots,B^{n_k-1}e'_k\}$;
then $B^{n_k}e'_k=\omega_k e'_k$ where $\omega_k$ is an $(\ell/n_k)$'th root of unity,
$V_k$ is a $B$-invariant subspace of ${\Bbb C}^n$,  and ${\cal E}_k$
is a basis of $V_k$, with which the matrix of $B$ restricted to $V_k$ is
\begin{equation}\label{monomial cycle}
   B|_{V_k}= \begin{pmatrix}0&\cdots&0&\omega_k\\ 1&\ddots&\ddots&0\\
    \vdots&\ddots&\ddots&\vdots\\ 0&\cdots&1&0\end{pmatrix}_{n_k\times n_k}.
\end{equation}

(2)~ ${\Bbb C}^n=V_1\oplus\cdots\oplus V_m$,
the union ${\cal E}={\cal E}_1\cup\cdots\cup{\cal E}_m$ is a basis
of ${\Bbb C}^n$ and, with respect to the basis ${\cal E}$, the matrix of $B$ is
\begin{equation}\label{B=oplus}
 B=B|_{V_1}\oplus\cdots\oplus B|_{V_m}.
\end{equation}
\end{lemma}

{\bf Proof.}~ (1).~
Since the length of $\Gamma_k$ is $n_k$,
$\Gamma_k=\big\{{\Bbb C}e'_k, B{\Bbb C}e'_k, \cdots,
 B^{n_k-1}{\Bbb C}e'_k\big\}$ and
$B^{n_k}{\Bbb C}e'_k={\Bbb C}e'_k$, hence there is an $\omega_k\in{\Bbb C}$
such that $B^{n_k}e'_k=\omega_k e'_k$. Then it is clear that $V_k$ is $B$-invariant,
${\cal E}_k$ is a basis of $V_k$ and Eqn (\ref{monomial cycle}) is
the matrix of $B|_{V_k}$. Since $B^\ell=I_{n\times n}$,
we have $(B|_{V_k})^{\ell}=I_{n_k\times n_k}$;
but by Eqn (\ref{monomial cycle}),
$(B|_{V_k})^{n_k}=\omega_k I_{n_k\times n_k}$;
so $\omega_k$ is an $(\ell/n_k)$'th root of unity.

(2).~
Applying (1) to all orbits $\Gamma_1,\cdots,\Gamma_m$, by Eqn (\ref{diag basis})
one can check the conclusions in~(2) easily.\qed

\begin{proposition}\label{SC}
Let notation be as in Lemma \ref{B block}.
Assume that the following condition is satisfied:

{\bf (SC)}~
For any $e_j$ and $B^i$, if $B^i{\Bbb C}e_j={\Bbb C}e_j$ then $B^ie_j=e_j$.

\noindent  Then the matrix group $\langle C,B\rangle$ generated by $C$ and $B$
is a permutation matrix group.
\end{proposition}

{\bf Proof.}~  We keep the notations in Lemma \ref{B block} and its proof.
We have seen that $B^{n_k}{\Bbb C}e'_k={\Bbb C}e'_k$;
by the condotion~(SC) we have $B^{n_k}e'_k=e'_k$, i.e. $\omega_k=1$ and
\begin{equation}\label{orbit cycle}
   B|_{V_k}= \begin{pmatrix}0&\cdots&0&1\\ 1&\ddots&\ddots&0\\
    \vdots&\ddots&\ddots&\vdots\\ 0&\cdots&1&0\end{pmatrix}_{n_k\times n_k};
\end{equation}
hence $B\sum_{h=0}^{n_k-1}B^he'_k=\sum_{h=0}^{n_k-1}BB^he'_k=
\sum_{h=0}^{n_k-1}B^he'_k$.

Now we set $f=\sum_{k=1}^{m}\sum_{h=0}^{n_k-1}B^he'_k$; then
$$Bf=B\sum_{k=1}^{m}\sum_{h=0}^{n_k-1}B^he'_k
  =\sum_{k=1}^{m}\sum_{h=0}^{n_k-1}BB^he'_k
  =\sum_{k=1}^{m}\sum_{h=0}^{n_k-1}B^he'_k=f.
$$
By Eqn (\ref{permutation basis}) the set of the vectors:
$$ f,~ Cf,~\cdots,~C^{n-1}f, $$
is a basis of ${\Bbb C}^n$; and with respect to this basis
$C$ is a cycle permutation matrix. Further, by Eqn (\ref{conjugation}) we have
$$
 BC^if=BC^iB^{-1}Bf=C^{r^{-1}i}f;
$$
that is, with respect to the basis  $f, Cf,\cdots,C^{n-1}f$,
the $B$ is also a permutation matrix.
In conclusion, the matrix group $\langle C,B\rangle$ generated by $C$ and $B$
is a permutation matrix group. \qed

\medskip
We'll quote a result of \cite{C07} repeatedly, so state it as a lemma:

\begin{lemma}\label{C07 P4.2}
{\rm(\cite[Proposition 4.2]{C07})}~
If $\langle C,B\rangle$ is an abelian  permutation-like matrix group where
$C$ is a maximal cycle, then $B\in\langle C\rangle$. \qed
\end{lemma}

We state some group-theoretic information as a remark for later quotations.
We say that an action of a group ${\cal G}$ on a set $\cal X$ is free
if the stabilizer of any $X\in{\cal X}$ in ${\cal G}$ is trivial;
at that case, ${\cal X}$ is partitioned in to ${\cal G}$-orbits
such that each orbit is a regular ${\cal G}$-set (i.e. equivalent to
the set ${\cal G}$ on which the group ${\cal G}$ acts by left translation).

\begin{remark}\label{p^2 group}\rm
Let $p$ be an odd prime.
\begin{itemize}
\item[(1)] Let ${\cal C}=\langle C\rangle$ be a cyclic group of order $p^2$.
It is easy to see that ${\cal C}^p=\langle C^p\rangle$ where
${\cal C}^p=\{X^p\,|\, X\in{\cal C}\}$;
and, mapping $X\in{\cal C}$ to $X^p\in{\cal C}^p$ is
a surjective homomorphism from ${\cal C}$ onto
${\cal C}^p$ with kernel  ${\cal C}^p$.

\item[(2)]
Let ${\cal G}$ be a finite group containing a normal
cyclic subgroup ${\cal C}=\langle C\rangle$ of order $p^2$
such that ${\cal C}$ is self-centralized
(i.e.  the centralizer $C_{\cal G}({\cal C})={\cal C}$).
Then ${\cal G}/{\cal C}$ is isomorphic to a subgroup of the automorphism
group of  ${\cal C}$, hence to a subgroup of the multiplicative group
${\Bbb Z}_{p^2}^*$ which is a cyclic group of order $p(p-1)$;
so there are an $r\in{\Bbb Z}_{p^2}^*$ and a $B\in{\cal G}$ such that
${\cal G}=\langle B,C\rangle$ and $C^B=C^{r}$
where $C^B=B^{-1}CB$ denotes the conjugate of $C$ by $B$;
and the order ${\rm ord}(r)=|{\cal G}/{\cal C}|$.

\item[(3)]
If ${\rm ord}(r)\big|(p-1)$, then ${\rm ord}(B)=|{\cal G}/{\cal C}|$
and the action by conjugation of $B$ on ${\cal C}$ is equivalent to
the action of $r$ by multiplication on ${\Bbb Z}_{p^2}$, the latter is
denoted by $\mu_r$, i.e. $\mu_r(a)=ra$ for all $a\in{\Bbb Z}_{p^2}$;
it is easy to check that the group $\langle B\rangle $ acts freely by conjugation on
the difference set ${\cal C}-\{1\}$
\big(this is a specific case of \cite[Corollary 4.35]{I}\big).

\item[(4)]
If  ${\rm ord}(r)=|{\cal G}/{\cal C}|=p$, then we can choose $B$ such that
$B^p=1$ and $r=ap+1$ with $0<a<p$; further, replacing $B$ by a suitable power of $B$,
we can get $r=p+1$; see \cite[\S8 Proposition 10]{AB} for details.

\item[(5)]
If $B^p=1$ and $r=p+1$ as in (4), then it is easy to check that:
\begin{itemize}
\item[(i)]
$\langle B\rangle$ centralizes the subgroup
${\cal C}^p=\langle C^p\rangle$ of ${\cal C}$,
and acts freely by conjugation on the difference set ${\cal C}-{\cal C}^p$;
in particular,  ${\cal C}-{\cal C}^p$ is partitioned into $p-1$ conjugacy classes
by $\langle B\rangle$, the length of every class is $p$.
\item[(ii)]
For any $X\in{\cal C}$ the product
$\prod_{j=0}^{p-1}(X)^{B^{j}}\in{\cal C}^p$;
the mapping $X$ to $\prod_{j=0}^{p-1}(X)^{B^{j}}$
is a surjective homomorphism from ${\cal C}$ onto ${\cal C}^p$;
hence the homomorphism induces a bijection
from the set of the conjugacy classes in ${\cal C}-{\cal C}^p$
on to the set ${\cal C}^p-\{1\}$.
\end{itemize}\end{itemize}
\end{remark}

\section{Proof of Theorem \ref{theorem}}

If $p=2$ then $p^2=4$ and the conclusion of Theorem \ref{theorem}
has been checked in \cite{C07}.

In the following, we always assume that $p$ is an odd prime and
${\cal G}$ is a permutation-like matrix group of dimension $p^2$
which contains a normal cyclic subgroup $\langle C\rangle$
generated by a maximal cycle $C$; and prove that
${\cal G}$ is a permutation matrix group.
If ${\cal G}$ is abelian, by Lemma \ref{C07 P4.2},
${\cal G}={\cal C}$ which is a permutation matrix group.
So we further assume that ${\cal G}$ is non-abelian.

\medskip
Let $\lambda$ be a primitive $p^2$'th root of unity,
By Eqn (\ref{diag basis}) there is a basis $e_0,e_1,\cdots,e_{p^2-1}$
of ${\Bbb C}^{p^2}$ such that
$$
C^ie_j=\lambda^{ij}e_j,\qquad i,j=0,1,\cdots,p^2-1.
$$
Let ${\cal C}=\langle C\rangle$ and $q=|{\cal G}/{\cal C}|$.
By Lemma \ref{C07 P4.2}, ${\cal C}$ is self-centralized in ${\cal G}$;
by Remark~\ref{p^2 group}(2),
there are a $B\in{\cal G}$, an $r\in{\Bbb Z}_{p^2}^*$ and integers $s,t$ such that
\begin{itemize}
\item
 $q=p^\delta s$ where $\delta=0$ or $1$, ~ $p-1=st$;
\item
${\cal G}=\langle B\rangle\cdot\langle C\rangle$ and the quotient
$\langle B\rangle\big/\langle B\rangle\cap\langle C\rangle\cong
\langle r\rangle\le{\Bbb Z}_{p^2}^*$;
\item
$B^{-1}CB=C^r$, i.e. the action by conjugation of
$\langle B\rangle\big/\langle B\rangle\cap\langle C\rangle$ on $\langle C\rangle$
is equivalent to the action of $\langle \mu_r\rangle$ on
the residual set ${\Bbb Z}_{p^2}$,
where $\mu_{r}(a)=ra$ for all $a\in{\Bbb Z}_{p^2}$.
\end{itemize}
The residual set ${\Bbb Z}_{p^2}$ is a disjoint union of two $\mu_{r}$-stable subsets:
\begin{equation}\label{Gamma_0}
{\Bbb Z}_{p^2}=\Gamma_0\cup{\Bbb Z}_{p^2}^*,\quad
 \mbox{where}~~
 \Gamma_0={\Bbb Z}_{p^2}-{\Bbb Z}_{p^2}^*=\{0,p,2p,\cdots,(p-1)p\};
\end{equation}
in fact, $\Gamma_0$ corresponds to the subgroup ${\cal C}^p$ of ${\cal C}$.
We prove the theorem in three cases.

\medskip
{\it Case 1.}~ $\delta=0$, i.e. $q=s\,|\,(p-1)$.~
By Remark \ref{p^2 group}(3),
$|\langle B\rangle|=s$ and
$\langle B\rangle\cong\langle r\rangle\le{\Bbb Z}_{p^2}^*$.
The group $\langle\mu_r\rangle$ fixes $0\in{\Bbb Z}_{p^2}$, and
acts freely on both ${\Bbb Z}_{p^2}^*$ and $\Gamma_0 -\{0\}$.
There are $t$ orbits of $\langle\mu_r\rangle$ on $\Gamma_0 -\{0\}$;
taking representatives $v_1,\cdots,v_t$ from the $t$ orbits,
we can write the orbits of $\langle \mu_r\rangle$ on $\Gamma_0$ as follows:
$$
\Gamma_{00}=\{0\}, \,\Gamma_{01}=\{v_1,rv_1,\cdots,r^{s-1}v_1\},\,\cdots,\,
\Gamma_{0t}=\{v_t,rv_t,\cdots,r^{s-1}v_t\}.
$$
There are $m$ orbits of $\langle\mu_r\rangle$ on ${\Bbb Z}_{p^2}^*$
where $m=\frac{p(p-1)}{s}=pt$;
taking representatives $w_1,\cdots,w_m$ from these $m$ orbits,
we have the orbits of $\langle \mu_r\rangle$ on ${\Bbb Z}_{p^2}^*$ as follows:
$$
\Gamma_{1}=\{w_1,rw_1,\cdots,r^{s-1}w_1\},\,\cdots,\,
\Gamma_{m}=\{w_m,rw_m,\cdots,r^{s-1}w_m\}.
$$
Accordingly, we apply Lemma \ref{B block} and its notation
to get the basis of ${\Bbb C}^{p^2}$:
$$
 {\cal E}={\cal E}_{00}\cup{\cal E}_{01}\cup\cdots\cup{\cal E}_{0t}\cup
     {\cal E}_{1}\cup\cdots\cup{\cal E}_{m},
$$
and write matrices with respect to this basis. So
$$
C=1\oplus\Big(\mathop{\oplus}_{i=0}^{s-1}\lambda^{r^{i}v_1}\Big)\oplus\cdots
   \oplus\Big(\mathop{\oplus}_{i=0}^{s-1}\lambda^{r^{i}v_t}\Big)
   \oplus\Big(\mathop{\oplus}_{i=0}^{s-1}\lambda^{r^{i}w_1}\Big)\oplus\cdots
   \oplus\Big(\mathop{\oplus}_{i=0}^{s-1}\lambda^{r^{i}w_m}\Big),
$$
and there is an $s$'th root $\varepsilon_0$ of unity such that
$$
B=\varepsilon_0\oplus\overbrace{P\oplus \cdots\oplus P}^{t+m},
$$
where $P$ is the cycle matrix of dimension $s$
(see Lemma \ref{B block} and Eqn (\ref{orbit cycle})):
$$  P=
   \begin{pmatrix}0&\cdots&0&1\\ 1&\ddots&\ddots&0\\
          \vdots&\ddots&\ddots&\vdots\\ 0&\cdots&1&0\end{pmatrix}_{s\times s}\,.$$
Then the characteristic polynomial of $B$ is
$$
{\rm char}_B(x)=(x-\varepsilon_0)(x^s-1)^{t+m};
$$
since $B$ is similar to a permutation matrix,
by Lemma \ref{permutable}, we have $\varepsilon_0=1$.
So the condition (SC) of Proposition \ref{SC} is satisfied:

{\bf (SC)}~ {\em For any $e_j$ and $B^i$, if $B^i{\Bbb C}e_j={\Bbb C}e_j$
then $B^ie_j=e_j$; }

\noindent
hence ${\cal G}$ is a permutation matrix group.

\medskip
{\it Case 2.} $s=1$, i.e. $q=p$. By Remark \ref{p^2 group}(4),
we can assume that $B^p=1$ and the conjugation of $B$ on ${\cal C}$
is equivalent to the action of $\mu_{p+1}$ on ${\Bbb Z}_{p^2}$,
where $\mu_{p+1}a=(p+1)a$ for $a\in{\Bbb Z}_{p^2}$;
further, $\mu_{p+1}$ centralizes the subset $\Gamma_0$ in Eqn~(\ref{Gamma_0}),
and  $\langle\mu_{p+1}\rangle$ partitions ${\Bbb Z}_{p^2}^*$
into $p-1$ $\langle\mu_{p+1}\rangle$-orbits, see Remark \ref{p^2 group}(5.i);
take representatives $u_1,\cdots,u_{p-1}$ from each $\langle\mu_{p+1}\rangle$-orbit,
the $p-1$ orbits can be written as:
\begin{eqnarray*}
  &&\Gamma_1=\{u_1, \,(p+1)u_1,\,\cdots,\,(p+1)^{p-1}u_1\},\\
   &&\Gamma_2=\{u_2,\,(p+1)u_2,\,\cdots,\,(p+1)^{(p-1)}u_2\},\\
   &&\cdots   \\
    &&\Gamma_{p-1}=\{u_{p-1},\,(p+1)u_{p-1},\,\cdots,\,(p+1)^{(p-1)}u_{p-1}\}.
\end{eqnarray*}
Accordingly, we apply Lemma \ref{B block} and its notation
to get the basis of ${\Bbb C}^{p^2}$:
$$
 {\cal E}={\cal E}_{0}\cup{\cal E}_{1}\cup\cdots\cup{\cal E}_{p-1},
$$
where ${\cal E}_k$ is corresponding to $\Gamma_k$ for $k=1,\cdots,p-1$
as above, while  ${\cal E}_0$ is corresponding to $\Gamma_0$;
and we write matrices with respect to this basis. So
$$
 C=D_0\oplus D_1\oplus\cdots\oplus D_{p-1},
$$
where
\begin{equation}\label{D_0}
D_0=1\oplus\lambda^p\oplus\cdots\oplus\lambda^{(p-1)p},
\end{equation}
$$
 D_i=\lambda^{u_i}\oplus\lambda^{(p+1)u_i}\oplus\cdots\oplus\lambda^{(p+1)^{p-1}u_i},
 \qquad i=1,\cdots,p-1;
$$
and
$$
B=B_0 \oplus \overbrace{P\oplus\cdots\oplus P}^{p-1},
$$
where
\begin{equation}\label{B_0}
B_0=\varepsilon_0\oplus \varepsilon_1\oplus\cdots\oplus \varepsilon_{p-1}
\end{equation}
with $\varepsilon_0,\varepsilon_1,\cdots,\varepsilon_{p-1}$
being $p$'th roots of unity and $P$ is the cycle matrix of dimension $p$
 (see Lemma \ref{B block} and Eqn (\ref{orbit cycle})):
$$ P=
   \begin{pmatrix}0&\cdots&0&1\\ 1&\ddots&\ddots&0\\
          \vdots&\ddots&\ddots&\vdots\\ 0&\cdots&1&0\end{pmatrix}_{p\times p}\,.$$

Any element of ${\cal G}$ has the form $C^kB^h$, $0\le k\le p^2-1$, $0\le h\le p-1$;
and
$$
  C^kB^h=D_0^k B_0^h\oplus D_1^{k}P^h\oplus\cdots\oplus D_{p-1}^{k}P^h.
$$
Obviously,
\begin{equation}\label{D0B0}
D_0^k B_0^h=\varepsilon_0^h\oplus\lambda^{pk}\varepsilon_1^h\oplus\cdots
  \oplus\lambda^{p(p-1)k}\varepsilon_{p-1}^h.
\end{equation}
It is easy to calculate the characteristic polynomials:
$$
{\rm char}_{D_i^kP^h}(x)=x^p-\lambda^{\sum_{j\in\Gamma_i}jk}
=x^p-\lambda^{u_ipmk},\qquad i=1,\cdots,p-1;
$$
where $pm=1+(p+1)+\cdots+(p+1)^{p-1}=\frac{(p+1)^p-1}{p}$, hence $m$
is an integer coprime to $p$.
If $k\not\equiv 0\pmod p$, by Remark \ref{p^2 group}(5.ii),
$\lambda^{u_ipmk}$ for $i=1,\cdots,p-1$ are just all primitive
$p$'th root of unity. Thus the characteristic polynomial of the matrix $C^k B^h$ is
$$
{\rm char}_{C^kB^h}(x)=\begin{cases}
 {\rm char}_{D_0^k B_0^h}(x)\cdot (x^p-1)^{p-1}, & k\equiv 0\pmod p;\\
 {\rm char}_{D_0^k B_0^h}(x)\cdot \Phi_{p^2}(x), & k\not\equiv 0\pmod p;
 \end{cases}
$$
where $\Phi_{p^2}(x)$ denotes the $p^2$'th cyclotomic polynomial.
Since $C^kB^h$ is similar to a permutation matrix,
by Lemma \ref{permutable}, for any $k,h$ we obtain that
\begin{equation}\label{D_0B_0}
{\rm char}_{D_0^k B_0^h}(x)=
\begin{cases} x^p-1~\mbox{or}~(x-1)^{p}, & k\equiv 0\pmod p;\\
                       x^p-1, & k\not\equiv 0\pmod p.\end{cases}
\end{equation}

By Eqn~(\ref{D_0}), we can view $D_0$ as a maximal cycle of dimension $p$;
by Eqn~(\ref{B_0}), the matrix $B_0$ of dimension $p$ commutes with $D_0$;
by Lemma \ref{permutable}, from Eqn (\ref{D_0B_0}) we see that
the abelian matrix group $\langle D_0,B_0\rangle$ of dimension $p$
is a permutation-like matrix group of dimension $p$;
so, by Lemma \ref{C07 P4.2},  we have an integer $0\le \ell\le p-1$ such that
$$
(\varepsilon_0,~\varepsilon_1,~\cdots,~\varepsilon_{p-1})
   =(1,~\lambda^{p\ell},~\cdots,~\lambda^{p(p-1)\ell}).
$$
Then, from Eqn (\ref{D0B0}) it is easy to calculate the characteristic polynomial
$$
{\rm char}_{D_0^kB_0^h}(x) =
   (x-1)(x-\lambda^{p(k+\ell h)})\cdots(x-\lambda^{p(p-1)(k+\ell h)}),
$$
so
\begin{equation}\label{D'_0}
{\rm char}_{D_0^kB_0^h}(x)
 =\begin{cases} (x-1)^{p}, & k+\ell h \equiv 0\pmod p;\\
                       x^p-1, & k+\ell h\not\equiv 0\pmod p. \end{cases}
\end{equation}
Suppose that $0<\ell\le p-1$; taking $k\not\equiv 0\pmod p$,
we have an $h$ such that $k+\ell h \equiv 0\pmod p$,
then, by Eqns (\ref{D_0B_0}) and (\ref{D'_0}) we have
$$
 x^p-1={\rm char}_{D_0^kB_0^h}(x)=(x-1)^p,
$$
which is impossible. Thus $\ell=0$, i.e.
$(\varepsilon_0,\varepsilon_1,\cdots,\varepsilon_{p-1})=(1,1,\cdots,1)$.

Summarizing the above, we obtain that
$$
 B=I_{p\times p}\oplus \overbrace{P\oplus\cdots\oplus P}^{p-1}.
$$
Similar to Case 1, the group $\langle B\rangle$ satisfies the condition (SC)
of Proposition \ref{SC}:

{\bf (SC)}~ {\em For any $e_j$ and $B^i$, if $B^i{\Bbb C}e_j={\Bbb C}e_j$
then $B^ie_j=e_j$. }

\noindent
Thus ${\cal G}$ is a permutation matrix group.

\medskip
{\em Case 3.}~ $q=ps$ and $s>1$.
First we show that
\begin{equation}\label{B^ps}
B^{ps}=C^{ap}\quad\mbox{with}~~ 0\le a< p,\quad\mbox{hence}~
{\rm ord}(B)=\begin{cases}ps, & a=0;\\ p^2s, & 0<a<p. \end{cases}
\end{equation}
\noindent For: otherwise $B^{ps}\in \langle C\rangle-\langle C^p\rangle$,
then $\langle C\rangle=\langle B^{ps}\rangle$, hence $B$ centralizes $\langle C\rangle$,
which contradicts to that ${\cal G}$ is non-abelian.

It is easy to see that the group $\langle\mu_{r}\rangle$ of order $q=ps$
acts freely on ${\Bbb Z}_{p^2}^*$, hence partitions ${\Bbb Z}_{p^2}^*$ in to
$t$ orbits of length $q$;
taking representatives $w_1,\cdots,w_t$ from these~$t$ orbits,
we have the orbits of $\langle \mu_r\rangle$ on ${\Bbb Z}_{p^2}^*$ as follows:
$$
\Gamma_{1}=\{w_1,rw_1,\cdots,r^{q-1}w_1\},\,\cdots,\,
\Gamma_{t}=\{w_t,rw_t,\cdots,r^{q-1}w_t\}.
$$

On the other hand,  since the group
$\langle\mu_{r^{s}}\rangle$ of order $p$ centralizes $\Gamma_0$,
the group $\langle\mu_{r}\rangle$ fixes $0$ and
partitions $\Gamma_0-\{0\}$ in to $t$ orbits of length $s$
(cf. Case 2); taking representatives $v_1,\cdots,v_t$ from these~$t$ orbits,
we have the orbits of $\langle \mu_r\rangle$ on $\Gamma_0$ as follows:
$$
\Gamma_{00}=\{0\}, \,\Gamma_{01}=\{v_1,rv_1,\cdots,r^{s-1}v_1\},\,\cdots,\,
\Gamma_{0t}=\{v_t,rv_t,\cdots,r^{s-1}v_t\}.
$$

According to the orbits $\Gamma_{00}$, $\Gamma_{01}$, $\cdots$,
$\Gamma_{0t}$, $\Gamma_{1}$, $\cdots$, $\Gamma_{t}$,
we apply Lemma \ref{B block} and its notation
to get the basis of ${\Bbb C}^{p^2}$:
$$
 {\cal E}={\cal E}_{00}\cup{\cal E}_{01}\cup\cdots\cup{\cal E}_{0t}\cup
     {\cal E}_{1}\cup\cdots\cup{\cal E}_{t},
$$
and write matrices with respect to this basis. Then
$$
C=1\oplus\Big(\mathop{\oplus}_{i=0}^{s-1}\lambda^{r^{i}v_1}\Big)\oplus\cdots
   \oplus\Big(\mathop{\oplus}_{i=0}^{s-1}\lambda^{r^{i}v_t}\Big)
   \oplus\Big(\mathop{\oplus}_{i=0}^{q-1}\lambda^{r^{i}w_1}\Big)\oplus\cdots
   \oplus\Big(\mathop{\oplus}_{i=0}^{q-1}\lambda^{r^{i}w_t}\Big);
$$
because $p|v_j$ for $j=1,\cdots,t$ (see Eqn (\ref{Gamma_0})),
$\lambda^{r^iv_jp}=1$ for $0\le i\le s-1$ and
$1\le j\le t$, so
$$
C^{ap}=I_{p\times p}\oplus
  \Big(\mathop{\oplus}_{i=0}^{q-1}\lambda^{r^{i}w_1ap}\Big)\oplus\cdots
   \oplus\Big(\mathop{\oplus}_{i=0}^{q-1}\lambda^{r^{i}w_tap}\Big).
$$
And, there are complexes $\varepsilon_0$ and
$\varepsilon_j$, $\omega_j$ for $j=1,\cdots,t$ such that
$$
B=\varepsilon_0\oplus P_1\oplus \cdots\oplus P_t
     \oplus Q_1\oplus \cdots\oplus Q_t,
$$
where $P_j$, $Q_j$ are as described in Eqn (\ref{monomial cycle})
(see Lemma \ref{B block}) :
$$
 P_j= \begin{pmatrix}0&\cdots&0&\varepsilon_j\\ 1&\ddots&\ddots&0\\
          \vdots&\ddots&\ddots&\vdots\\ 0&\cdots&1&0\end{pmatrix}_{s\times s}.
          \qquad
 Q_j= \begin{pmatrix}0&\cdots&0&\omega_j\\ 1&\ddots&\ddots&0\\
          \vdots&\ddots&\ddots&\vdots\\ 0&\cdots&1&0\end{pmatrix}_{q\times q}\,;
$$
hence
$$
 B^{ps}=\varepsilon_0^{ps}\oplus
    \varepsilon_1^p I_{s\times s}\oplus\cdots\oplus\varepsilon_t^p I_{s\times s}
     \oplus\omega_1 I_{q\times q}\oplus\cdots\oplus\omega_t I_{q\times q} .
$$
Since $B^{ps}=C^{ap}$ \big(see Eqn (\ref{B^ps})\big), the collection of
\begin{equation}\label{omega}
\overbrace{\omega_1,~\cdots,~\omega_1}^q,~~\cdots,~~
\overbrace{\omega_t,~\cdots,~\omega_t}^q
\end{equation}
is coincide with the collection of
\begin{equation}\label{ap}
 \lambda^{w_1ap}, \lambda^{rw_1ap}, \cdots, \lambda^{r^{q-1}w_1ap},~\cdots,~
  \lambda^{w_tap}, \lambda^{rw_tap}, \cdots, \lambda^{r^{q-1}w_tap}.
\end{equation}
Note that
$$
\big\{ w_1,  rw_1, \cdots, r^{q-1}w_1,~\cdots,~
  w_t, rw_t, \cdots, r^{q-1}w_t\big\}={\Bbb Z}_{p^2}^*.
$$
If $0<a<p$ (cf. Eqn (\ref{B^ps})),
then the collection (\ref{ap}) is just all primitive $p$'th roots of unity
with multiplicity $p$ for each one, see Remark \ref{p^2 group}(1);
on the other hand, the number of the elements
appeared in the collection (\ref{omega}) is at most  $t$;
but $ts=p-1$ and $s>1$,
so $t$ is less than the number of primitive $p$'th roots;
that is a contradiction to
the coincidence of the collections (\ref{omega}) and (\ref{ap}).

In conclusion, $a=0$ and $B^{ps}=1$.

Since $p$ and $s$ are coprime, we have
$$\langle B\rangle=\langle B^p\rangle\times\langle B^s\rangle,\quad
    |\langle B^p\rangle|=s,\quad |\langle B^s\rangle|=p.
$$
We have considered $\langle C,B^p\rangle$ and $\langle C,B^s\rangle$
in Case 1 and Case 2 respectively, and have concluded that
in both cases the condition (SC) in Proposition \ref{SC} is satisfied.
Assume that $B^i{\Bbb C}e_j={\Bbb C}e_j$;
taking integers $k,h$ such that $ph+sk=1$, we have $B^i=B^{phi}B^{ski}$,
$$
  B^{phi}\in\langle B^p\rangle,~  B^{phi}{\Bbb C}e_j={\Bbb C}e_j~~{\rm and}~~
  B^{ski}\in\langle B^s\rangle,~ B^{ski}{\Bbb C}e_j={\Bbb C}e_j;
$$
by the conclusions in Case 1 and in Case 2,
$B^{phi}e_j=e_j$ and $B^{ski}e_j=e_j$; hence $B^i e_j=e_j$.
Thus, by Proposition \ref{SC},
 ${\cal G}$ is a permutation matrix group.

The proof of Theorem \ref{theorem} is completed.

\section*{Acknowledgements}
The research of the authors is supported by NSFC
with grant numbers 11171194 and 11271005.

\end{document}